\documentclass{amsart}
 
\usepackage{amsmath,amssymb,amsthm}
\usepackage{enumitem}
 
\newtheorem{theorem}{Theorem}
\newtheorem{corollary}[theorem]{Corollary}
\newtheorem{proposition}[theorem]{Proposition}
\newtheorem{lemma}[theorem]{Lemma}
\theoremstyle{remark}
\newtheorem{remark}[theorem]{Remark}
\newtheorem{example}[theorem]{Example}
 
\DeclareMathOperator{\Span}{span}
\DeclareMathOperator{\Tr}{Tr}
\DeclareMathOperator{\rank}{rank}
\DeclareMathOperator{\End}{End}
\DeclareMathOperator{\row}{row}
 
\newcommand{\C}{\mathbb C}
\newcommand{\norm}[1]{\lVert #1\rVert}
 
\title[Irreducible representations as local maxima]{Irreducible Representations as Local Maxima of a Matrix-Coefficient Functional}
\author{Pin-Chi Hung}
\address{Department of Mathematics, Soochow University, Taipei City, Taiwan}
\email{pinchihung1111@gmail.com}
 
\author{Ming-Hsuan Kang}
\address{Department of Applied Mathematics, National Yang Ming Chiao Tung University, Hsinchu, Taiwan}
\email{kmsming@gmail.com}
 
 
\begin{document}
 
\begin{abstract}
Let $V$ be a finite-dimensional unitary representation of a finite group $G$ and, for unit vectors
$v$, let $\Phi(v)$ be the mean of $|\langle v,gv\rangle|^2$ over $g\in G$. We show that $v$ is a local
maximum of $\Phi$ if and only if the subrepresentation generated by $v$ is irreducible. Thus local
maximization of $\Phi$ extracts irreducible subrepresentations of $V$ directly, with no character or
isotypic projection specified in advance. The proof is elementary: Schur averaging reduces it to a
singular-value calculation.
The values and dimensions of the connected components of the local-maximum set determine the degrees
and multiplicities of the irreducible constituents, and for $H\le G$ the intersections of the $G$- and $H$-local-maximal loci
encode branching multiplicities. We also determine the critical landscape and give a normalized
gradient iteration whose generic limit generates an irreducible subrepresentation.
\end{abstract}
\maketitle
 
\section{Introduction}
 
Let $G$ be a finite group and let $\rho$ be a representation of $G$ on a finite-dimensional complex
vector space $V$, unitary with respect to a $G$-invariant inner product. Characters determine the
isotypic components of $V$ through central idempotents, but when a multiplicity exceeds one they do
not single out an irreducible subrepresentation. Character-free methods do: one may average a
Hermitian operator into the commutant and diagonalize it, as in Dixon's method \cite{Dixon1970}, or
use more refined algorithms \cite{BabaiRonyai1990}. We ask whether irreducible subrepresentations can
instead be located by the local geometry of a single scalar function on $V$, built only from the
matrix coefficients of one vector. Consider
\begin{equation}\label{eq:Phi}
    \Phi(v)=\frac1{|G|}\sum_{g\in G}\bigl|\langle v,\rho(g)v\rangle\bigr|^2,
    \qquad \norm{v}=1.
\end{equation}
Our main result is
\begin{equation}\label{eq:main-intro}
 v\text{ is a local maximum of }\Phi
 \quad\Longleftrightarrow\quad
 U_v:=\Span\{\rho(g)v:g\in G\}\text{ is irreducible}.
\end{equation}
Thus the local maxima of $\Phi$ are exactly the generators of irreducible subrepresentations. The
point is the variational description rather than a new
decomposition algorithm: no irreducible type is prescribed in advance, and the connected components
of the local-maximum set correspond to the irreducible types occurring in $V$, with their values and
dimensions recovering the degrees and multiplicities.
 
Our second result concerns restriction. Let $H\le G$. The same vector space carries both the $G$- and
$H$-functionals, and Theorem~\ref{thm:branching} shows that the projectivized local-maximal components
attached to irreducible representations $\pi$ of $G$ and $\sigma$ of $H$ meet in a product of three
projective spaces, whose middle factor has dimension one less than the multiplicity of $\sigma$ in the
restriction of $\pi$. Thus the incidence and dimensions of these intersections recover the branching
multiplicities; in the regular representation they recover the full restriction matrix. Frobenius
reciprocity and Gelfand pairs acquire direct interpretations in the same geometry.
 
The proof is elementary. After identifying the isotypic components of $V$ with matrix spaces, the
squared singular values of the component matrices of a unit vector form a probability vector. Schur
averaging expresses $\Phi$ as a weighted sum of squares of its entries. If two entries are positive,
transferring a small amount of mass between them increases $\Phi$ in one direction. Hence a local
maximum has exactly one nonzero entry, which is precisely the irreducible case. The same computation determines the full critical set and shows
that $\Phi$ is Morse--Bott (Section~\ref{sec:critical}), and it yields a normalized gradient iteration,
using only inner products of orbit vectors, that reaches a local maximum from a generic starting vector
(Section~\ref{sec:iter}).
 
The averaging step is classical. In quantum information it is the group twirl
\cite{BartlettRudolphSpekkens2007}, whose relation to asymmetry is well known
\cite{MarvianSpekkens2014}. In frame-theoretic terms, the critical points are the vectors with
tight orbits, classified by Vale and Waldron \cite[Theorem~6.18]{ValeWaldron2005}; in that setting tight
frames usually appear as minimizers of frame potentials \cite{BenedettoFickus2003,JohnsonOkoudjou2008},
whereas here the emphasis is on local maxima. This note grew out of \cite{HsiehKang2026}, where an
irreducible character is prescribed in advance and $\Phi$ is maximized on the corresponding isotypic
component; the upper bound for $\Phi$ on that component and its equality case are proved there. What
is new here is the global statement: with no character prescribed, the local maxima of $\Phi$ on all
of $V$ are exactly the generators of irreducible subrepresentations. The description of the
local-maximum set, the branching theorem, and the Morse--Bott classification of the critical points
have no counterpart in \cite{HsiehKang2026}.

\section{The averaging operator}\label{sec:avg}
 
Throughout, inner products are linear in the second variable, and derivatives on $V$ are
taken with respect to the real inner product $\operatorname{Re}\langle\cdot,\cdot\rangle$. When $V$ is
reducible, different $G$-invariant inner products give different functionals; the results below hold
for each of them.
 
\paragraph{Matrix model.}
Let $(\rho_\alpha)_\alpha$ be unitary representatives of the irreducible representations of $G$
occurring in $V$, of degrees $d_\alpha$ and multiplicities $m_\alpha$. The isotypic components of
$V$ are mutually orthogonal, and the $\alpha$-component is an orthogonal sum of $m_\alpha$
invariant subspaces, each unitarily $G$-isomorphic to $\C^{d_\alpha}$ \cite{Serre1977}. Using these
isomorphisms we identify
\[
 V=\bigoplus_\alpha \C^{d_\alpha\times m_\alpha},\qquad
 \rho(g)(X_\alpha)_\alpha=(\rho_\alpha(g)X_\alpha)_\alpha,\qquad
 \langle X,Y\rangle=\sum_\alpha\Tr(X_\alpha^*Y_\alpha),
\]
the $j$-th column of $X_\alpha$ being the component in the $j$-th copy. For
$v=(X_\alpha)_\alpha$ put
\[
 U_v=\Span\{\rho(g)v:g\in G\}.
\]
Lemma~\ref{lem:key} below shows that the averaging operator acts on the $\alpha$-component through
the positive matrix $X_\alpha^*X_\alpha$ on the multiplicity space. Consequently, its eigenvalues, the
value of $\Phi$, and the local geometry of $\Phi$ are all governed by the singular values
$s_{\alpha,j}$ of the matrices $X_\alpha$. We put
\[
 t_{\alpha,j}=s_{\alpha,j}^2,\qquad r_\alpha=\rank X_\alpha.
\]
If $v$ is a unit vector, then $\sum_{\alpha,j}t_{\alpha,j}=1$.
 
\paragraph{The averaging operator.}
For $v\in V$, let $vv^*$ denote the rank-one operator $u\mapsto v\langle v,u\rangle$ and put
\begin{equation}\label{eq:Kv}
 K_v=\frac1{|G|}\sum_{g\in G}\rho(g)vv^*\rho(g)^*.
\end{equation}
We equip $\End(V)$ with the Hilbert--Schmidt inner product
$\langle A,B\rangle_{\mathrm{HS}}=\Tr(A^*B)$, and write $\norm{A}_F^2=\Tr(A^*A)$ for the Frobenius
norm of a rectangular matrix. Then $K_v$ is the orthogonal projection of $vv^*$ onto the commutant
$\End_G(V)=\{A\in\End(V):A\rho(g)=\rho(g)A\text{ for all }g\in G\}$.
 
\begin{lemma}[Key lemma]\label{lem:key}
Let $v=(X_\alpha)_\alpha$ be a unit vector.
\begin{enumerate}[label=\textup{(\alph*)}]
\item $K_v$ acts on $\C^{d_\alpha\times m_\alpha}$ by
\[
 Y\longmapsto d_\alpha^{-1}YX_\alpha^*X_\alpha.
\]
\item $K_v$ is positive semidefinite, commutes with $\rho(G)$, has trace $1$, and has range $U_v$.
Its nonzero eigenvalues are $t_{\alpha,j}/d_\alpha$ $(1\le j\le r_\alpha)$, each with multiplicity
$d_\alpha$. More precisely,
\[
 U_v=\bigoplus_\alpha U_{v,\alpha},\qquad
 U_{v,\alpha}
 =\{Y\in\C^{d_\alpha\times m_\alpha}:\row(Y)\subseteq\row(X_\alpha)\}
 \cong \rho_\alpha^{\oplus r_\alpha},
\]
and therefore
\[
 \dim U_v=\sum_\alpha d_\alpha r_\alpha.
\]
\item
\[
 \Phi(v)=\langle v,K_vv\rangle=\norm{K_v}_{\mathrm{HS}}^2
 =\sum_{\alpha,j}\frac{t_{\alpha,j}^2}{d_\alpha}.
\]
\item The gradient of the quartic form $\langle v,K_vv\rangle$ is $4K_vv$.
\end{enumerate}
\end{lemma}
 
\begin{proof}
The $\alpha\beta$ block of the average \eqref{eq:Kv} is a $G$-intertwiner between the
$\beta$- and $\alpha$-isotypic components, so Schur's lemma kills the blocks with
$\alpha\ne\beta$. If $x_c,x_e$ are columns of $X_\alpha$, then
\[
 \frac1{|G|}\sum_{g\in G}
 \rho_\alpha(g)x_cx_e^*\rho_\alpha(g)^*
 =\frac{\langle x_e,x_c\rangle}{d_\alpha}I_{d_\alpha}.
\]
Thus $K_v$ acts on the $\alpha$-component through the multiplicity space, by right multiplication
by $d_\alpha^{-1}X_\alpha^*X_\alpha$:
\[
 K_vY=d_\alpha^{-1}YX_\alpha^*X_\alpha
 \qquad(Y\in\C^{d_\alpha\times m_\alpha}),
\]
which proves (a).
 
For (b), if $w_g=\rho(g)v$, then
\[
 \langle u,K_vu\rangle=\frac1{|G|}\sum_g|\langle w_g,u\rangle|^2,
\]
so $K_v\ge0$ and $\ker K_v=U_v^\perp$; hence $\operatorname{ran}K_v=U_v$. Also
$\Tr K_v=1$, and reindexing the average shows that $K_v$ commutes with $\rho(G)$. Part (a) now gives the
stated eigenvalues and
\[
 U_{v,\alpha}=\{Y:\row(Y)\subseteq\row(X_\alpha)\}
 \cong\rho_\alpha^{\oplus r_\alpha}.
\]
Part (c) follows at once from (a):
\[
 \langle v,K_vv\rangle
 =\sum_\alpha d_\alpha^{-1}\Tr((X_\alpha^*X_\alpha)^2)
 =\sum_{\alpha,j}\frac{t_{\alpha,j}^2}{d_\alpha},
\]
and the same expression is $\norm{K_v}_{\mathrm{HS}}^2$ by (b). Finally,
$\nabla\Tr((X^*X)^2)=4XX^*X$, which together with (a) gives (d).
\end{proof}

\begin{remark}[Extension to compact groups]
Lemma~1 remains valid when $G$ is a compact Hausdorff group
and $\rho:G\to U(V)$ is a continuous unitary representation
on a finite-dimensional complex inner product space $V$.
Let $\mu$ be normalized Haar measure on $G$, and replace
the finite averages by
\[
\Phi(v)
  = \int_G |\langle v,\rho(g)v\rangle|^2\,d\mu(g),
\qquad
K_v
  = \int_G \rho(g)vv^*\rho(g)^*\,d\mu(g).
\]
These integrals are well defined by continuity and compactness.
Haar averaging is again the Hilbert--Schmidt orthogonal
projection onto $\operatorname{End}_G(V)$, and Schur
orthogonality gives the same isotypic and singular-value
formulas as in Lemma~1. In particular, for $\|v\|=1$,
\[
\Phi(v)=\sum_{\alpha,j}\frac{t_{\alpha,j}^2}{d_\alpha},
\qquad
\sum_{\alpha,j}t_{\alpha,j}=1.
\]
Consequently, the proof of the local-maximum characterization
below also applies: a unit vector $v$ is a local maximum
of $\Phi$ if and only if its cyclic subrepresentation $U_v$
is irreducible, in which case $\Phi(v)=1/\dim U_v$.
\end{remark}
 
\section{Local maxima}\label{sec:max}
 
Lemma~\ref{lem:key} suggests working in singular coordinates. For each $\alpha$, choose a
singular value decomposition
\begin{equation}\label{eq:svd}
 X_\alpha
 =\sum_{j=1}^{r_\alpha}s_{\alpha,j}u_{\alpha,j}w_{\alpha,j}^*
 =\sum_{j=1}^{r_\alpha}\sqrt{t_{\alpha,j}}\,\xi_{\alpha,j},
 \qquad
 \xi_{\alpha,j}:=u_{\alpha,j}w_{\alpha,j}^*.
\end{equation}
The rank-one matrices $\xi_{\alpha,j}$ are orthonormal in $\C^{d_\alpha\times m_\alpha}$, and
different isotypic components are orthogonal. Hence
\begin{equation}\label{eq:v-singular}
 v=\sum_{\alpha,j}\sqrt{t_{\alpha,j}}\,\xi_{\alpha,j},
 \qquad
 \sum_{\alpha,j}t_{\alpha,j}=1.
\end{equation}
By Lemma~\ref{lem:key}(c),
\begin{equation}\label{eq:Phi-singular}
 \Phi(v)=\sum_{\alpha,j}\frac{t_{\alpha,j}^2}{d_\alpha}.
\end{equation}
With the singular directions $\xi_{\alpha,j}$ held fixed, $\Phi$ is therefore a weighted quadratic
function on the simplex of masses $(t_{\alpha,j})$. On the other hand, by Lemma~\ref{lem:key}(b),
$U_v$ is irreducible exactly when $\sum_\alpha r_\alpha=1$, that is, when a single mass
$t_{\alpha,j}$ is nonzero (and hence equal to $1$).
 
For a single prescribed isotypic component, the bound $\Phi\le 1/d$ and its equality case were
established in \cite{HsiehKang2026}. The following theorem is the global counterpart: local
maximality alone forces a vector into one isotypic component and one irreducible copy.
 
\begin{theorem}\label{thm:localmax}
A unit vector $v$ is a local maximum of $\Phi$ if and only if $U_v$ is irreducible. In that case
\[
 \Phi(v)=\frac1{\dim U_v}.
\]
\end{theorem}
 
\begin{proof}
Suppose two of the masses in \eqref{eq:v-singular} are positive; call them $t_1,t_2$, with
corresponding singular directions $\xi_1,\xi_2$ and degrees $d_1,d_2$ (possibly $d_1=d_2$). Write
\[
 v=\sqrt{t_1}\,\xi_1+\sqrt{t_2}\,\xi_2+v_0,
 \qquad
 v_0\perp \xi_1,\xi_2.
\]
For $|\varepsilon|<\min(t_1,t_2)$, we transfer mass between the two directions by setting
\[
 v(\varepsilon)=\sqrt{t_1+\varepsilon}\,\xi_1+\sqrt{t_2-\varepsilon}\,\xi_2+v_0.
\]
Since the directions in \eqref{eq:v-singular} are orthonormal, $\norm{v(\varepsilon)}=1$ and $v(0)=v$, and
all other masses are unchanged. By \eqref{eq:Phi-singular},
\[
 \Phi(v(\varepsilon))-\Phi(v)
 =2\varepsilon\left(\frac{t_1}{d_1}-\frac{t_2}{d_2}\right)
 +\varepsilon^2\left(\frac1{d_1}+\frac1{d_2}\right).
\]
Choosing the sign of $\varepsilon\ne0$ so that the linear term is nonnegative, we obtain
$\Phi(v(\varepsilon))>\Phi(v)$ for arbitrarily small $|\varepsilon|$, because the quadratic coefficient is positive.
Hence $v$ is not a local maximum. It follows that a local maximum has exactly one nonzero mass, so
$U_v$ is irreducible.
 
Conversely, suppose that $U_v$ is irreducible, so that $v$ has a single nonzero mass
$t_{\alpha,1}=1$, lying in a component of degree $d=d_\alpha$. For a unit vector $w$ near $v$, let
$1-\delta$ be the largest eigenvalue of $X_\alpha(w)^*X_\alpha(w)$; by continuity, $\delta\ge0$ is
small. The remaining masses $t_{\beta,k}(w)$ sum to $\delta$, and $d_\beta\ge1$ for all $\beta$.
Hence
\[
 \Phi(w)\le \frac{(1-\delta)^2}{d}+\delta^2\le \frac1d=\Phi(v)
 \qquad\text{whenever }0\le\delta\le \frac2{d+1}.
\]
Therefore $v$ is a local maximum.
\end{proof}
 
\begin{corollary}\label{cor:components}
The local-maximum set of $\Phi$ is the disjoint union, over the irreducible constituents $\rho_\alpha$
of $V$, of the sets $\mathcal M_\alpha$ of unit rank-one matrices in
$\C^{d_\alpha\times m_\alpha}$. Moreover,
\[
 \mathcal M_\alpha\cong
 \bigl(S^{2d_\alpha-1}\times S^{2m_\alpha-1}\bigr)/S^1,
 \qquad
 \dim_{\mathbb R}\mathcal M_\alpha=2d_\alpha+2m_\alpha-3,
 \qquad
 \Phi|_{\mathcal M_\alpha}=\frac1{d_\alpha}.
\]
In particular, the connected components of the local-maximum set are in bijection with the irreducible
constituents of $V$ and determine their degrees and multiplicities.
\end{corollary}
 
\begin{proof}
The first statement restates Theorem~\ref{thm:localmax}. A unit rank-one matrix can be written $uw^*$ with
$u\in S^{2d_\alpha-1}$ and $w\in S^{2m_\alpha-1}$, uniquely up to multiplying $u$ and $w$ by the
same phase, and this $S^1$-action is free. Since the sets $\mathcal M_\alpha$ are compact,
connected and pairwise disjoint, they are the connected components.
\end{proof}
 
The condition in Theorem~\ref{thm:localmax} does not involve the inner product. Hence the image of
the local-maximum set in $\mathbb P(V)$, namely the set of lines generating irreducible subrepresentations,
depends only on the $G$-module $V$ and not on the choice of invariant inner product; the same holds
for the values $1/d_\alpha$. The remaining critical points do depend on this choice.
 
\subsection{Subgroups and branching}\label{subsec:branching}
 
Let $H\le G$. The same representation $V$ carries two functionals, according to whether the average is
taken over $G$ or over $H$. To distinguish them, write
\[
 K_v^L=\frac1{|L|}\sum_{\ell\in L}\rho(\ell)vv^*\rho(\ell)^*,
 \qquad
 \Phi_L(v)=\norm{K_v^L}_{\mathrm{HS}}^2,
 \qquad L=G,H,
\]
so that $K_v^L=\mathcal E_L(vv^*)$ for the averaging map
$\mathcal E_L(A)=|L|^{-1}\sum_{\ell\in L}\rho(\ell)A\rho(\ell)^*$ on $\End(V)$.
The scalar potential is monotone under restriction, and the defect is the squared distance between the
two averaged operators.
 
The map $\mathcal E_L$ is the orthogonal projection of $\End(V)$ onto $\End_L(V)$. Since
$\End_G(V)\subseteq\End_H(V)$, we have $\mathcal E_G\mathcal E_H=\mathcal E_G$. Thus $K_v^G$ is the orthogonal projection of $K_v^H$
onto $\End_G(V)$, and Pythagoras gives
\begin{equation}\label{eq:subgroup-refinement}
 \Phi_H(v)-\Phi_G(v)
 =\norm{K_v^H-K_v^G}_{\mathrm{HS}}^2\ge0.
\end{equation}
At a common local maximum of types $\pi\in\widehat G$ and $\sigma\in\widehat H$, the defect equals
$1/\dim\sigma-1/\dim\pi$ by Theorem~\ref{thm:localmax}.
 
The local-maximal loci encode the branching data of the $G$-types occurring in $V$. Write the $G$-isotypic decomposition as
\begin{equation}\label{eq:G-isotypic-branching}
 V\cong\bigoplus_{\pi\in\widehat G}V_\pi\otimes M_\pi,
 \qquad m_\pi=\dim M_\pi,
\end{equation}
and, for each $\pi$, write its restriction to $H$ as
\begin{equation}\label{eq:restriction-branching}
 \operatorname{Res}_H^G V_\pi
 \cong\bigoplus_{\sigma\in\widehat H}W_\sigma\otimes B_{\sigma\pi},
 \qquad b_{\sigma\pi}=\dim B_{\sigma\pi}.
\end{equation}
Hence, as an $H$-module,
\[
 V\cong\bigoplus_{\sigma\in\widehat H}W_\sigma\otimes M_\sigma^H,
 \qquad
 M_\sigma^H=\bigoplus_{\pi\in\widehat G}B_{\sigma\pi}\otimes M_\pi.
\]
Since $\Phi_L(e^{i\theta}v)=\Phi_L(v)$, the set of local maxima of $\Phi_L$ descends to a subset $\mathbb P\mathcal M^L(V)$ of
$\mathbb P(V)$. Its components corresponding to $\pi\in\widehat G$ and $\sigma\in\widehat H$ are denoted
$\mathbb P\mathcal M_\pi^G(V)$ and $\mathbb P\mathcal M_\sigma^H(V)$; we omit $(V)$ when the
ambient module is clear. If an $H$-type $\sigma$ does not
occur in $V$, we use the convention $\mathbb P\mathcal M_\sigma^H=\varnothing$.
 
The two local-maximal loci need not be nested. The following criterion describes when an
entire $G$-local-maximal component survives after restriction to $H$.

\begin{proposition}[Survival under restriction]\label{prop:component-survival}
Let $\pi\in\widehat G$ be a type occurring in $V$. Then
\[
 \mathbb P\mathcal M_\pi^G\subseteq \mathbb P\mathcal M^H(V)
\]
if and only if either $\operatorname{Res}_H^G V_\pi$ is irreducible or $H$ acts on $V_\pi$ through a
single character.
\end{proposition}

\begin{proof}
A point of $\mathbb P\mathcal M_\pi^G$ is represented by $u\otimes m$, with $0\ne u\in V_\pi$
and $0\ne m\in M_\pi$. Its $H$-span is the $H$-span of $u$ tensored with $m$, so by
Theorem~\ref{thm:localmax} it remains $H$-locally maximal exactly when the $H$-span of $u$ is
irreducible. If two distinct $H$-types occur in \eqref{eq:restriction-branching}, a vector with nonzero
components in both generates a reducible $H$-module. Hence all nonzero $u$ can work only when
$\operatorname{Res}_H^G V_\pi\cong W_\sigma\otimes B_{\sigma\pi}$ for a single $\sigma$. In this case,
by Lemma~\ref{lem:key}(b) applied to $H$, the number of copies of $W_\sigma$ generated by $u$ is the
matrix rank of $u\in W_\sigma\otimes B_{\sigma\pi}$. Thus every nonzero $u$ generates an irreducible
$H$-module exactly when
\[
 \min\{\dim W_\sigma,\dim B_{\sigma\pi}\}=1.
\]
If $\dim B_{\sigma\pi}=1$, the restriction is irreducible; if $\dim W_\sigma=1$, the group $H$
acts on $V_\pi$ through the character $\sigma$.
\end{proof}

Thus only smaller pieces of a $G$-local-maximal component usually survive, while an $H$-local
maximum may involve several $G$-types. The next theorem describes the surviving pieces exactly.
 
\begin{theorem}[Branching from common local maxima]\label{thm:branching}
For $\pi\in\widehat G$ occurring in $V$ and $\sigma\in\widehat H$,
\[
 \mathbb P\mathcal M_\pi^G\cap\mathbb P\mathcal M_\sigma^H=\varnothing
 \qquad\Longleftrightarrow\qquad
 b_{\sigma\pi}=0.
\]
If $b_{\sigma\pi}>0$, then
\begin{equation}\label{eq:branching-intersection}
 \mathbb P\mathcal M_\pi^G\cap\mathbb P\mathcal M_\sigma^H
 \cong
 \mathbb P(W_\sigma)\times
 \mathbb P(B_{\sigma\pi})\times
 \mathbb P(M_\pi).
\end{equation}
In particular,
\begin{equation}\label{eq:branching-dimension}
 \dim_{\C}\bigl(\mathbb P\mathcal M_\pi^G\cap\mathbb P\mathcal M_\sigma^H\bigr)
 =d_\sigma+b_{\sigma\pi}+m_\pi-3,
 \qquad d_\sigma=\dim W_\sigma.
\end{equation}
Thus the incidence of the two local-maximal loci recovers which $H$-types occur in
$\operatorname{Res}_H^G V_\pi$. Together with the degree and multiplicity already recorded by the
individual components, the intersection dimensions recover the branching multiplicities
$b_{\sigma\pi}$.
\end{theorem}
 
\begin{proof}
By Corollary~\ref{cor:components}, a point of $\mathbb P\mathcal M_\pi^G$ is represented by a simple
tensor $u\otimes m$ with $u\in V_\pi$ and $m\in M_\pi$. Under \eqref{eq:restriction-branching}, its
$\sigma$-isotypic part lies in
\[
 W_\sigma\otimes B_{\sigma\pi}\otimes M_\pi.
\]
For the same vector to lie in $\mathbb P\mathcal M_\sigma^H$, it must have $H$-type $\sigma$ and be
rank one for the decomposition $W_\sigma\otimes M_\sigma^H$. Inside the fixed $\pi$-summand, the two
rank-one conditions are the two flattenings
\[
 (W_\sigma\otimes B_{\sigma\pi})\otimes M_\pi,
 \qquad
 W_\sigma\otimes(B_{\sigma\pi}\otimes M_\pi).
\]
A nonzero tensor has rank one in both flattenings exactly when it is fully decomposable. Indeed, if
$t=a\otimes m=w\otimes c$ with $a\in W_\sigma\otimes B_{\sigma\pi}$ and
$c\in B_{\sigma\pi}\otimes M_\pi$, choose $\lambda\in M_\pi^*$ with $\lambda(m)=1$. Contracting the
third factor gives
\[
 a=w\otimes (\operatorname{id}\otimes\lambda)(c),
\]
so $a=w\otimes b$ and hence $t=w\otimes b\otimes m$; the converse is immediate. This gives
\eqref{eq:branching-intersection}; the isomorphism is induced by the natural Segre embedding
$([w],[b],[m])\mapsto[w\otimes b\otimes m]$, so the intersection is a smooth irreducible projective
variety. It is empty precisely when
$B_{\sigma\pi}=0$. Formula \eqref{eq:branching-dimension} follows by adding the three projective
dimensions. Finally, $d_\sigma$ is read from the value
$\Phi_H=1/d_\sigma$ on $\mathbb P\mathcal M_\sigma^H$, while $m_\pi$ is already determined by the
$G$-local-maximal component, so \eqref{eq:branching-dimension} determines $b_{\sigma\pi}$.
\end{proof}

The most intrinsic form occurs in the left regular representation $V=\C[G]$, where the multiplicity space $M_\pi$ has dimension $d_\pi$.
Hence Theorem~\ref{thm:branching} gives, whenever $b_{\sigma\pi}>0$,
\begin{equation}\label{eq:regular-branching}
 \mathbb P\mathcal M_\pi^G\cap\mathbb P\mathcal M_\sigma^H
 \cong
 \mathbb P^{d_\sigma-1}\times
 \mathbb P^{b_{\sigma\pi}-1}\times
 \mathbb P^{d_\pi-1}.
\end{equation}
Thus the two local-maximal loci on the same space $\C[G]$ recover all restriction multiplicities,
with rows and columns indexed by the corresponding $H$- and $G$-components.
 
\begin{example}\label{ex:S4-S2}
Take $G=S_4$ and $H=\langle(12)\rangle\cong S_2$. The standard representation
$V_{(31)}$ has degree $3$, and its character takes the value $1$ on a transposition. Hence
\[
 \operatorname{Res}_{S_2}^{S_4}V_{(31)}\cong 2\mathbf1\oplus\operatorname{sgn}.
\]
In $\C[S_4]$, the common local-maximal locus for $V_{(31)}$ and the trivial $S_2$-type is therefore
\[
 \mathbb P^0\times\mathbb P^1\times\mathbb P^2,
\]
of complex dimension $3$. The middle factor records the branching multiplicity $2$.
\end{example}
 
\subsection{Frobenius reciprocity and Gelfand pairs}\label{subsec:gelfand}
 
For $\sigma\in\widehat H$, Frobenius reciprocity gives
$[\operatorname{Ind}_H^G W_\sigma:V_\pi]=b_{\sigma\pi}$. Hence, by Corollary~\ref{cor:components}, the
$\pi$-component of the local-maximum set of $\operatorname{Ind}_H^G W_\sigma$ has multiplicity factor
$\mathbb P(\C^{b_{\sigma\pi}})$ when $b_{\sigma\pi}>0$, matching the middle factor in
Theorem~\ref{thm:branching}. For $\sigma=\mathbf1_H$, the induced module is
$\Omega=\C[G/H]$. Define
\[
 \mathcal S(G,H)
 =\mathbb P\mathcal M^G(\Omega)\cap\mathbb P(\Omega^H),
 \qquad
 \mathcal S_\pi(G,H)
 =\mathbb P\mathcal M_\pi^G(\Omega)\cap\mathbb P(\Omega^H),
\]
so that $\mathcal S(G,H)$ is the disjoint union of the sets $\mathcal S_\pi(G,H)$ over the
$\pi\in\widehat G$ occurring in $\Omega$.
 
\begin{corollary}[Gelfand pairs]\label{cor:gelfand}
Write
\[
 \Omega\cong\bigoplus_{\pi\in\widehat G}V_\pi\otimes M_\pi,
 \qquad m_\pi=\dim M_\pi.
\]
For every occurring $\pi$,
\begin{equation}\label{eq:gelfand-component}
 \mathcal S_\pi(G,H)
 \cong
 \mathbb P(V_\pi^H)\times\mathbb P(M_\pi)
 \cong
 \mathbb P^{m_\pi-1}\times\mathbb P^{m_\pi-1}.
\end{equation}
Consequently,
\[
 (G,H)\text{ is a Gelfand pair}
 \quad\Longleftrightarrow\quad
 \mathcal S(G,H)\text{ is finite}.
\]
When these conditions hold, each irreducible constituent of $\C[G/H]$ contributes one point to
$\mathcal S(G,H)$. If $v_\pi$ is a unit $H$-fixed representative of that point, then
\[
 \varphi_\pi(g)=\langle v_\pi,\rho(g)v_\pi\rangle
\]
is the corresponding normalized zonal spherical function.
\end{corollary}
 
\begin{proof}
Frobenius reciprocity gives $m_\pi=\dim V_\pi^H$. On the $\pi$-isotypic component
$V_\pi\otimes M_\pi$, the group $H$ acts only on the first factor. Corollary~\ref{cor:components}
therefore gives
\[
 \mathcal S_\pi(G,H)
 =\mathbb P\mathcal M_\pi^G(\Omega)\cap\mathbb P(V_\pi^H\otimes M_\pi)
 \cong\mathbb P(V_\pi^H)\times\mathbb P(M_\pi),
\]
which proves \eqref{eq:gelfand-component}. This is zero-dimensional exactly when $m_\pi=1$. Hence
$\mathcal S(G,H)$ is finite exactly when $\C[G/H]$ is multiplicity-free. In that case the unique
$H$-fixed line in each irreducible summand gives the normalized zonal spherical function $\varphi_\pi$.
\end{proof}
 
\section{Critical points and the Morse--Bott structure}\label{sec:critical}
 
We begin with a lower bound that also identifies the critical values.
 
\begin{theorem}\label{thm:critical}
For every unit vector $v$,
\begin{equation}\label{eq:critical-lower}
 \Phi(v)\ge\frac1{\dim U_v}.
\end{equation}
Equality holds exactly at the critical points of $\Phi$. More precisely, the following are equivalent:
\begin{enumerate}[label=\textup{(\roman*)}]
\item $v$ is a critical point of $\Phi$ on the unit sphere;
\item $K_vv=\Phi(v)v$;
\item $K_v=P_{U_v}/\dim U_v$;
\item all nonzero numbers $t_{\alpha,j}/d_\alpha$ are equal.
\end{enumerate}
In this case $\Phi(v)=1/D$ and every nonzero singular value of $X_\alpha$ equals
$\sqrt{d_\alpha/D}$, where $D=\dim U_v=\sum_\alpha d_\alpha r_\alpha$.
\end{theorem}
 
\begin{proof}
By Lemma~\ref{lem:key}, $K_v$ has $D=\dim U_v$ nonzero eigenvalues, counted with multiplicity,
summing to $1$, while $\Phi(v)=\norm{K_v}_{\mathrm{HS}}^2$ is the sum of their squares. This gives
\eqref{eq:critical-lower}, with equality exactly when $K_v=P_{U_v}/D$.
 
On the unit sphere, Lemma~\ref{lem:key}(d) shows that $v$ is critical exactly when
$K_vv=\lambda v$; pairing with $v$ gives $\lambda=\Phi(v)$. If this holds, commutation with $\rho(G)$
shows that $K_v$ is scalar on $U_v=\operatorname{ran}K_v$ and zero on $U_v^\perp$, so its trace
forces $K_v=P_{U_v}/D$. The converse is immediate, and Lemma~\ref{lem:key}(b) identifies this with
(iv). The final assertions follow from the common eigenvalue $1/D$.
\end{proof}
 
By Lemma~\ref{lem:key}(b), a critical point with rank data $(r_\alpha)$ generates
$U_v\cong\bigoplus_\alpha\rho_\alpha^{\oplus r_\alpha}$. It therefore generates an irreducible
subrepresentation exactly when $\sum_\alpha r_\alpha=1$, that is, exactly when it is a local maximum.
 
\begin{remark}[Frame-theoretic interpretation]\label{rem:frames}
The orbit $(\rho(g)v)_{g\in G}$ has frame operator $|G|K_v$ and frame potential $|G|^2\Phi(v)$.
Thus Theorem~\ref{thm:critical} says that the critical points are precisely the vectors whose orbit is a
tight frame for its span $U_v$. The equivalence of \textup{(iii)} and \textup{(iv)} is then the
classification of tight group orbits of Vale and Waldron \cite[Theorem~6.18]{ValeWaldron2005}, applied
to $U_v$; see also \cite{BenedettoFickus2003,Han2007}.
\end{remark}
 
We now describe the critical set as a union of manifolds.
For $0\le r\le\min(d,m)$, let $\mathcal P_r(d,m)$ be the set of rank-$r$ partial isometries in
$\C^{d\times m}$, that is, matrices $X$ for which $X^*X$ is an orthogonal projection of rank $r$
(so $\mathcal P_0=\{0\}$). It is a single orbit of the compact connected group $U(d)\times U(m)$
acting by $X\mapsto UXV^*$, and hence a compact connected submanifold. By
Theorem~\ref{thm:critical}, the critical set of $\Phi$ is the disjoint union of the manifolds
\[
 \mathcal C_{\mathbf r}
 =\prod_\alpha \sqrt{\frac{d_\alpha}{D_{\mathbf r}}}\,\mathcal P_{r_\alpha}(d_\alpha,m_\alpha),
 \qquad
 D_{\mathbf r}=\sum_\alpha d_\alpha r_\alpha,
\]
over rank data $\mathbf r=(r_\alpha)$ with $0\le r_\alpha\le\min(d_\alpha,m_\alpha)$, not all zero.
The critical value on $\mathcal C_{\mathbf r}$ is $1/D_{\mathbf r}$.
 
Recall that a smooth function is \emph{Morse--Bott} if its critical set is a union of smooth
submanifolds and, at every critical point, the kernel of the Hessian equals the tangent space of the
critical submanifold through that point. We write $n_+$ and $n_-$ for the numbers of positive and
negative Hessian directions normal to that submanifold.
 
\begin{theorem}\label{thm:morse}
$\Phi$ is a Morse--Bott function on the unit sphere of $V$. Along $\mathcal C_{\mathbf r}$, the Hessian has
\[
 n_+=\sum_\alpha r_\alpha^2-1
 \qquad\text{positive and}\qquad
 n_-=2\sum_\alpha(d_\alpha-r_\alpha)(m_\alpha-r_\alpha)
\]
negative normal directions. In particular, the local maxima are the points with
$\sum_\alpha r_\alpha^2=1$ (recovering Theorem~\ref{thm:localmax}), and the local minima are the
points with $r_\alpha=\min(d_\alpha,m_\alpha)$ for all $\alpha$.
\end{theorem}
 
\begin{proof}
By \eqref{eq:Phi-singular}, $\Phi$ depends only on singular values, so $\Phi$, the unit sphere and
$\mathcal C_{\mathbf r}$ are all invariant under $X_\alpha\mapsto U_\alpha X_\alpha V_\alpha^*$ with
$U_\alpha\in U(d_\alpha)$ and $V_\alpha\in U(m_\alpha)$. We may therefore assume
\[
 X_\alpha=s_\alpha
 \begin{pmatrix}I_{r_\alpha}&0\\0&0\end{pmatrix},
 \qquad s_\alpha^2=\frac{d_\alpha}{D},\qquad D=D_{\mathbf r}.
\]
Write a tangent vector in the same block form as
\[
 Z_\alpha=\begin{pmatrix}A_\alpha&B_\alpha\\C_\alpha&E_\alpha\end{pmatrix}.
\]
Since $\nabla\Phi(v)=4K_vv=(4/D)v$, the constrained second derivative on the unit sphere is
$D^2\Phi(v)[Z,Z]-(4/D)\norm{Z}^2$. For $f(X)=\Tr((XX^*)^2)$ one has
\[
 D^2f(X)[Z,Z]
 =2\norm{XZ^*+ZX^*}_F^2+4\Tr(XX^*ZZ^*).
\]
In block form, for $X=s\bigl(\begin{smallmatrix}I&0\\0&0\end{smallmatrix}\bigr)$,
\[
 \norm{XZ^*+ZX^*}_F^2=s^2\bigl(\norm{A+A^*}_F^2+2\norm{C}_F^2\bigr),
 \qquad
 \Tr(XX^*ZZ^*)=s^2\bigl(\norm{A}_F^2+\norm{B}_F^2\bigr).
\]
At our point, $s_\alpha^2/d_\alpha=1/D$, so the terms in $\norm{A_\alpha}_F^2$, $\norm{B_\alpha}_F^2$
and $\norm{C_\alpha}_F^2$ cancel against $(4/D)\norm{Z}^2$, and the Hessian on the sphere becomes
\[
 \operatorname{Hess}\Phi(v)[Z,Z]
 =\frac2D\sum_\alpha\norm{A_\alpha+A_\alpha^*}_F^2
 -\frac4D\sum_\alpha\norm{E_\alpha}_F^2.
\]
Write $A_\alpha=S_\alpha+Q_\alpha$, with $S_\alpha^*=S_\alpha$ and
$Q_\alpha^*=-Q_\alpha$. The sphere tangent condition becomes
\[
 \sum_\alpha s_\alpha\Tr S_\alpha=0,
\]
and the Hessian is
\[
 \frac8D\sum_\alpha\norm{S_\alpha}_F^2
 -\frac4D\sum_\alpha\norm{E_\alpha}_F^2.
\]
The Hermitian blocks $S_\alpha$ contribute $\sum_\alpha r_\alpha^2-1$ positive directions after
the tangent constraint, while the blocks $E_\alpha$ contribute the stated negative directions. The kernel
is given by $S_\alpha=E_\alpha=0$. These are exactly the tangent directions to the
$U(d_\alpha)\times U(m_\alpha)$ orbit: differentiating
$X_\alpha\mapsto e^{tL_\alpha}X_\alpha e^{-tR_\alpha}$ gives
\[
 s_\alpha\begin{pmatrix}
 L_{\alpha,11}-R_{\alpha,11}&-R_{\alpha,12}\\
 L_{\alpha,21}&0
 \end{pmatrix},
\]
which supplies all skew-Hermitian upper-left and off-diagonal directions. Hence the Hessian kernel is
$T_v\mathcal C_{\mathbf r}$, proving the Morse--Bott claim and the final assertions.
\end{proof}
 
\begin{corollary}[Landscape]\label{cor:landscape}
Put $d_{\min}=\min_\alpha d_\alpha$ and
$D_{\max}=\sum_\alpha d_\alpha\min(d_\alpha,m_\alpha)$. Then
\[
 \min\Phi=\frac1{D_{\max}},\qquad \max\Phi=\frac1{d_{\min}}.
\]
Every local minimum is global; the global maxima are exactly the local maxima in constituents of degree
$d_{\min}$; and every remaining critical manifold is a saddle.
\end{corollary}
 
\begin{proof}
The minimum statement follows from Theorems~\ref{thm:critical} and~\ref{thm:morse}, since
$\dim U_v\le D_{\max}$, with equality exactly at the maximal rank data. For the maximum,
\eqref{eq:Phi-singular} gives
\[
 \Phi(v)\le d_{\min}^{-1}\sum_{\alpha,j}t_{\alpha,j}^2\le d_{\min}^{-1},
\]
with equality exactly when one mass is nonzero in a constituent of degree $d_{\min}$. The saddle
statement is immediate from the signs in Theorem~\ref{thm:morse}.
\end{proof}
 
Local maxima need not be global: on the component of $\rho_\alpha$ their value is $1/d_\alpha$.

\section{The normalized gradient iteration}\label{sec:iter}
 
Theorem~\ref{thm:localmax} identifies irreducible cyclic subrepresentations with local maxima.
The same averaging operator gives a simple iteration toward these maxima. Define
\[
 T(v)=\frac{K_vv}{\norm{K_vv}}.
\]
This is well defined on the unit sphere, since $\langle v,K_vv\rangle=\Phi(v)>0$, and by
Lemma~\ref{lem:key}(d) it is the normalized ambient-gradient map.
 
\begin{lemma}[Commutant power identity]\label{lem:power}
For every unit vector $v$,
\begin{equation}\label{eq:K-cube}
 K_{T(v)}=\frac{K_v^3}{\Tr(K_v^3)}.
\end{equation}
Consequently, if $v^{(0)}=v$, $v^{(n+1)}=T(v^{(n)})$, and $K_n=K_{v^{(n)}}$, then
\begin{equation}\label{eq:K-closed}
 K_n=\frac{K_0^{3^n}}{\Tr(K_0^{3^n})}.
\end{equation}
\end{lemma}
 
\begin{proof}
Let $\mathcal E=\mathcal E_G$, so that $K_v=\mathcal E(vv^*)$. For
$B,C\in\End_G(V)$, $\mathcal E(BAC)=B\mathcal E(A)C$, and therefore
\[
 K_{T(v)}
 =\frac{\mathcal E(K_vvv^*K_v)}{\norm{K_vv}^2}
 =\frac{K_v^3}{\norm{K_vv}^2}.
\]
Since $\mathcal E$ preserves trace,
$\norm{K_vv}^2=\Tr(\mathcal E(vv^*K_v^2))=\Tr(K_v^3)$, proving
\eqref{eq:K-cube}; iteration gives \eqref{eq:K-closed}.
\end{proof}
 
\begin{corollary}\label{cor:iteration}
Suppose the largest scaled singular value $s_{\alpha,j}/\sqrt{d_\alpha}$ of the isotypic matrices of
$v$ occurs at a unique index $(\alpha,j)$. Then $T^n(v)$ converges to a local maximum generating an
irreducible subrepresentation of type $\rho_\alpha$.
\end{corollary}
 
\begin{proof}
By Lemma~\ref{lem:key}, the eigenvalues of $K_v$ are the numbers
$s_{\beta,k}^2/d_\beta$. Formula \eqref{eq:K-closed} raises their ratios to the power $3^n$.
Moreover, iterating $v^{(n+1)}\propto K_n v^{(n)}$ gives
\[
 v^{(n)}\propto K_0^{(3^n-1)/2}v
 \propto
 \sum_{\beta,k}\sqrt{d_\beta}
 \left(\frac{s_{\beta,k}}{\sqrt{d_\beta}}\right)^{3^n}\xi_{\beta,k}.
\]
Thus the unique dominant term survives, giving a rank-one limit; Theorem~\ref{thm:localmax} identifies
it as a local maximum of type $\rho_\alpha$.
\end{proof}
 
The hypothesis of Corollary~\ref{cor:iteration} holds outside a proper real-algebraic subset of the
unit sphere, hence on an open dense set of full measure. Indeed, let $\Gamma_\alpha$ be the smaller of
$X_\alpha^*X_\alpha$ and $X_\alpha X_\alpha^*$, of size $\min(d_\alpha,m_\alpha)$, and let $c_\alpha$
be the characteristic polynomial of $d_\alpha^{-1}\Gamma_\alpha$. Its roots are the squares of the
scaled singular values of $X_\alpha$, zeros included. Put
\[
 P(v)=\prod_\alpha\operatorname{disc}(c_\alpha)
 \prod_{\alpha<\beta}\operatorname{Res}(c_\alpha,c_\beta),
\]
a homogeneous real polynomial in the coordinates of $v$. If $P(v)\ne0$, all these roots are distinct,
so the largest scaled singular value occurs at a unique index. Moreover, $P$ is not identically zero:
it is nonzero when every $X_\alpha$ is rectangular diagonal and the numbers
$|(X_\alpha)_{jj}|^2/d_\alpha$, $1\le j\le\min(d_\alpha,m_\alpha)$, are nonzero and pairwise distinct
over all $(\alpha,j)$. Hence the iteration converges to an irreducible local maximum from every unit
vector with $P(v)\ne0$. Every irreducible type occurring in $V$ is the limit type for a nonempty open
set of starting vectors, namely those near a rank-one vector of that type.

For the left regular representation on the group algebra $V=\C[G]$, with the inner product for which
$G$ is an orthonormal basis, $K_x$ is right multiplication by $|G|^{-1}x^*x$, where
$(\sum_ga_gg)^*=\sum_g\overline{a_g}g^{-1}$, so the
iteration becomes simply
\[
 T(x)=\frac{xx^*x}{\norm{xx^*x}}.
\]
Splitting a representation by averaging and then diagonalizing an element of the commutant goes back
to Dixon \cite{Dixon1970}; the cubing identity isolates the dominant eigenspace without
diagonalization.

\section*{Acknowledgements}
The authors used OpenAI ChatGPT (GPT-5.6 Sol; OpenAI) for exploratory mathematical discussion, proof
checking and language revision during the development of this work.  The authors take full responsibility
for the mathematical content and its verification.

\end{document}